# Proof of a conjecture of N. Konno for the 1D contact process

## J. van den Berg[1],[*], O. Häggström[2] and J. Kahn[3],[†]

*CWI and VUA, Chalmers University of Technology and Rutgers University*

**Abstract:** Consider the one-dimensional contact process. About ten years ago, N. Konno stated the conjecture that, for all positive integers $n, m$, the upper invariant measure has the following property: Conditioned on the event that $O$ is infected, the events {All sites $-n, \ldots, -1$ are healthy} and {All sites $1, \ldots, m$ are healthy} are negatively correlated.

We prove (a stronger version of) this conjecture, and explain that in some sense it is a dual version of the planar case of one of our results in [2].

## 1. Introduction and statement of the main result

Consider the contact process on $\mathbf{Z}$ with infection rates $\lambda(x, y)$, $x, y \in \mathbf{Z}$, $|x - y| = 1$, and recovery rates $\delta_x, x \in \mathbf{Z}$. This model can, somewhat informally, be described as follows: Each site $x \in \mathbf{Z}$ has, at each time $t \geq 0$, a value $\eta_t(x) \in \{0, 1\}$. Usually, 1 is interpreted as 'infected' (or 'ill') and 0 as 'healthy'. When a site $y$ is ill, it infects each neighbour $x$ at rate $\lambda(x, y)$. In other words, at time $t$ each healthy site $x$ becomes infected at rate $\lambda(x, x+1)\eta_t(x+1) + \lambda(x, x-1)\eta_t(x-1)$. Further, when a site $x$ is ill, it recovers (becomes healthy) at rate $\delta_x$.

We assume that the above-mentioned rates are bounded. In fact, the most commonly studied case is where all recovery rates are constant, say 1, and all infection rates are equal to some value $\lambda > 0$. See [6] and [7] for background and further references.

Let $\nu_t$ be the law of $(\eta_t(x), x \in \mathbf{Z})$. It is well known that if at time 0 all sites are infected, $\nu_t$ converges, as $t \to \infty$, to a distribution called the upper invariant measure. We denote this limit distribution by $\nu$.

About ten years ago N. Konno proposed the following conjecture (see [4], Conjecture 4.5.2 or [5], Conjecture 2.3.2):

**Conjecture 1.** *Let $\nu$ be the upper invariant measure for the 1D contact process with infection rate $\lambda$ and recovery rate 1. Let $n, m$ be positive integers. Then*

$$\nu(\eta(x) = 0, \, x = -n, \ldots, -1, 1, \ldots, m \,|\, \eta(0) = 1)$$
$$\leq \nu(\eta(x) = 0, \, x = -n, \ldots, -1 \,|\, \eta(0) = 1)$$
$$\times \nu(\eta(x) = 0, \, x = 1, \ldots, m \,|\, \eta(0) = 1). \tag{1}$$

---

[*]Part of JvdB's research is financially supported by BRICKS project AFM2.2.

[†]Supported in part by NSF grant DMS0200856.

[1]CWI (Department PNA), Kruislaan 413, 1098 SJ Amsterdam, The Netherlands, e-mail: J.van.den.Berg@cwi.nl

[2]Chalmers University of Technology, S-412 96 Göteborg, Sweden, e-mail: olleh@math.chalmers.se

[3]Department of Mathematics, Rutgers University, Piscataway NJ 07059, USA, e-mail: jkahn@math.rutgers.edu

Before we state our stronger version, we give some notation and terminology. A finite collection, say $X_1, \ldots, X_n$, of $0 - 1$ valued random variables is said to be positively associated if, for all functions $f, g : \{0, 1\}^n \to \mathbf{R}$ that are both increasing or both decreasing,

$$E(f(X_1, \ldots, X_n) \, g(X_1, \ldots, X_n)) \geq E(f(X_1, \ldots, X_n)) E(g(X_1, \ldots, X_n)). \quad (2)$$

Equivalently, if $f$ is increasing and $g$ decreasing (or vice versa),

$$E(f(X_1, \ldots, X_n) \, g(X_1, \ldots, X_n)) \leq E(f(X_1, \ldots, X_n)) E(g(X_1, \ldots, X_n)). \quad (3)$$

Further, a countable collection of random variables is said to be positively associated, if every finite subcollection is positively associated.

We are now ready to state our main result, a stronger version of Konno's conjecture.

**Theorem 2.** *Let $\eta_t(x), x \in \mathbf{Z}, t \geq 0$ be the 1D contact process with deterministic initial configuration, and bounded infection and recovery rates. For each $t$ we have that, conditioned on the event that $\eta_t(0) = 1$, the collection of random variables $\{1 - \eta_t(x) \, : \, x < 0\} \cup \{\eta_t(x) \, : \, x > 0\}$ is positively associated.*

**Remark.** This theorem easily implies Conjecture 1: Start the contact process with all sites infected. Let $t > 0$ and let $n, m$ be positive integers. Let $A$ be the event $\{\eta_t(x) = 0, -n \leq x \leq -1\}$, and $B$ the event $\{\eta_t(x) = 0, 1 \leq x \leq m\}$. The indicator function of $A$ is an increasing function of the tuple $(1 - \eta_t(x), -n \leq x \leq -1)$, and the indicator function of $B$ is a decreasing function of the tuple $(\eta_t(x), 1 \leq x \leq m)$. Hence, by Theorem 2, conditioned on the event $\{\eta_t(0) = 1\}$, the events $A$ and $B$ are negatively correlated. This holds for each $t$. The conjecture follows by letting $t \to \infty$.

## 2. Slight extension of an earlier inequality

In this section we present a slight extension of a result in [2]. In Section 3 we will prove (for certain graphs) a dual version of this extension.

Let $G$ be a finite, or countably infinite, mixed graph. The word 'mixed' means that we allow that some of the edges are oriented and others non-oriented. A non-oriented edge between vertices $x$ and $y$ is denoted by $\{x, y\}$, and an oriented edge from $x$ to $y$ by $(x, y)$.

Let $V = V(G)$ and $E = E(G)$ be the vertex and edge sets of $G$. Let $p(e), e \in E$, be values in $[0, 1]$. Consider the percolation model on $G$ where each edge $e$, independently of the others, is open with probability $p_e$ and closed with probability $1 - p_e$.

When we speak of a path in $G$, we assume that it respects the orientation of its edges. As usual, an open path is a path every edge of which is open.

For $S, T \subset V(G)$, the event that there is an open path from (some vertex in) $S$ to (some vertex in) $T$ will be denoted by $\{S \to T\}$, the complement of this event by $\{S \nrightarrow T\}$, and the indicators of these events by $I_{S \to T}$ and $I_{S \nrightarrow T}$ respectively.

**Theorem 3.** *Let $S$ and $T$ be disjoint subsets of $V(G)$. Let, for each edge $e$, $X_e$ and $Y_e$ be the indicators of the events $\{e$ belongs to an open path beginning in $S\}$ and $\{e$ belongs to an open path ending in $T\}$ respectively. Then, conditioned on the event $\{S \nrightarrow T\}$, the collection $\{X_e : e \in E(G)\} \cup \{1 - Y_e : e \in E(G)\}$ is positively associated.*

This is an oriented generalization of Theorem 1.5 of [2] and follows from a straightforward modification of the arguments in that paper. (See Section 3 of [2] for other generalizations).

## 3. A planar dual version of Theorem 3

When the graph $G$ in Section 2 is embeddable in the plane, one can obtain a 'planar dual version' of Theorem 3. This (for the case of finite $G$, which for our purpose is sufficient) is Theorem 4 below. In Section 4 we will apply Theorem 4 to a special graph, which can be regarded as a discrete-time version of the ususal space–time diagram of the contact process. This application will yield Theorem 2.

**Theorem 4.** *Let $G = (V, E)$ be a finite, planar, mixed graph. Let $C$ be a subset of $E$ which, when one disregards edge orientations, forms a face-bounding cycle in some planar embedding of $G$. Let $u_1, \ldots, u_k, a_1, \ldots, a_m, w_1, \ldots, w_l$ and $b_1, \ldots, b_n$ denote (in some cyclic order) the vertices of $C$. Let $U \subset \{u_1, \ldots, u_k\}$ and $W \subset \{w_1, \ldots, w_l\}$. Consider bond percolation on $G$ with parameters $p_e, e \in E$. Then, conditioned on the event $\{U \to W\}$, the collection of random variables*

$$\{I_{U \to a_i} : 1 \leq i \leq m\} \cup \{I_{U \not\to b_j} : 1 \leq j \leq n\}$$

*is positively associated.*

*Proof.* We assume, without loss of generality that $C$ bounds the outer face in the given embedding of $G$, that the vertices of $C$ are given above in clockwise order, and that $k = \ell = 1$. (To justify the last assertion, add vertices $u, w$ outside $C$ and undirected edges joining $u$ to the $u_i$'s in $U$ and $w$ to the $w_i$'s in $W$, and let these new edges be open with probability 1.) So we may simply take $U = \{u\}$ and $W = \{w\}$, and condition on $\{u \to w\}$.

If, for some $x, y \in V(G)$, we have $(x, y) \in E(G)$ but $(y, x) \notin E(G)$, we can just add $(y, x)$ to $E(G)$ and take $p_{(y,x)} = 0$ without essentially changing anything. So (again, w.l.o.g.) we assume that $(x, y) \in E(G)$ iff $(y, x) \in E(G)$.

Finally, if $\{x, y\}$ is an undirected edge, which is open with probability $p$, we replace this edge by two directed edges which are independently open with probability $p$. It is well-known and easy to check that this does not change the distribution of the collection $(I_{u \to v}, v \in V(G))$, and hence it does not affect the assertion of Theorem 4. Therefore, we assume w.l.o.g. that all edges of $G$ are directed.

Next, for convenience, we slightly vary the usual definition of the undirected graph, $\overline{G}$, underlying $G$. The graph $\overline{G}$ has the same vertices as $G$. All edges of $\overline{G}$ are undirected, and $\{x, y\}$ is an edge of $\overline{G}$ iff $(x, y)$ (and hence, by one of the assumptions above, also $(y, x)$) is an edge of $G$.

It is clear that the alternative form of Theorem 4 obtained by replacing $G$ by $\overline{G}$ in the second line, is equivalent to the original form. From now on we will refer to that alternative form.

The following conventions also turn out to be convenient: we consider a "drawing" (not, strictly speaking, *embedding*) of $G$ which coincides with the given embedding of $\overline{G}$, in the sense that $(x, y), (y, x) \in E(G)$ are both drawn as orientations of the curve representing the corresponding edge in the embedding of $\overline{G}$. These conventions will also apply to the dual-like graph $H$ defined below. The dual $e^*$ of an edge $e$ will always be oriented to cross $e$ from left to right (as these sides are understood when one follows the direction of $e$).

For $x, y \in V(C)$ we use $[x, y]$ for the set of edges of $G$ whose underlying edges (in $\overline{G}$) belong to the path obtained by following $C$ clockwise from $x$ to $y$.

We form a graph $H$, a variant of the planar dual of $G$, as follows. Start with vertices corresponding to the bounded faces of (our drawing of) $G$, joining them by dual edges as usual (oriented according to the preceding convention, and again taking $(x, y)^*$ and $(y, x)^*$ to be represented by the same curve). Then, for each $e \in E(G)$ with underlying edge belonging to $C$, add a dual edge $e^*$ joining the vertex of $H$ corresponding to the inner face containing $e$ in its boundary to a new vertex $s_e$ lying in the outer face. The $s_e$'s are distinct except that $s_{(x,y)} = s_{(y,x)}$. (To avoid introducing unwanted crossings, take $s_e$ to be drawn just outside $e$.)

For $x \in V(C) \setminus \{u\}$, let $S_x = \{s_e : e \in [u, x]\}$, and $T_x = \{s_e : e \in [x, u]\}$, and set $S_w = S$, $T_w = T$.

We couple percolation on $H$ with that on $G$ in the natural way, by declaring an edge $e^*$ of $H$ to be open (closed) if the corresponding edge $e$ of $G$ is closed (open).

Let $V^*$ and $E^*$ denote the vertex set and the edge set of $H$ respectively. For connection events in $H$ we use similar notation as for $G$, with the symbol '$*$' added to indicate that we consider the dual. For instance, if $s$ and $t$ are vertices of $H$, $s \overset{*}{\to} t$ denotes the event that there is an open path in $H$ from $s$ to $t$, and $s \overset{*}{\not\to} t$ the complement of that event.

We will apply Theorem 3 to the graph $H$.

*Observations*

(i) For each $x \in V(C) \setminus \{u\}$, $u \to x$ iff $S_x \overset{*}{\not\to} T_x$. (This is an analog of standard duality properties of planar percolation). In particular, $u \to w$ iff $S \overset{*}{\not\to} T$.

(ii) Let, for each edge $e^*$ of $H$, $X_{e^*}$ and $Y_{e^*}$ be as in Theorem 3 (i.e. $X_{e^*}$ is the indicator of the event that $e^*$ belongs to an an open path in $H$ beginning in $S$, and $Y_{e^*}$ is the indicator of the event that $e^*$ belongs to an open path in $H$ ending in $T$). By observation (i), for each $i$, $1 \le i \le m$, $I_{u \to a_i}$ is a decreasing function of the collection $(X_{e^*} : e^* \in E^*)$, and for each $j$, $1 \le j \le n$, $I_{u \to b_j}$ is a decreasing function of the collection $(Y_{e^*} : e^* \in E^*)$.

Now let $f$ and $g$ be increasing functions of the collection $\{I_{u \to a_i} : 1 \le i \le m\} \cup \{I_{u \to b_j} : 1 \le j \le n\}$. By observation (ii), $f$ and $g$ are decreasing functions of the collection $\{X_{e^*} : e^* \in E^*\} \cup \{1 - Y_{e^*} : e^* \in E^*\}$. We get:

$$
\begin{aligned}
E(fg \,|\, u \to w) &= E(fg \,|\, S \overset{*}{\not\to} T) \qquad\qquad (4) \\
&\ge E(f \,|\, S \overset{*}{\not\to} T) \, E(g \,|\, S \overset{*}{\not\to} T) \\
&= E(f \,|\, u \to w) \, E(g \,|\, u \to w),
\end{aligned}
$$

where the equalities follow from observation (i), and the inequality follows from Theorem 3. This completes the proof of Theorem 4. $\qquad \square$

### 3.1. *An alternative proof of Theorem 4*

We think, and we believe this is also part of Mike's philosophy, that a problem is best understood by approaching it in several ways. This subsection gives a sketch of a self-contained proof of Theorem 4. Instead of using explicit *results* from [2], it uses *ideas* similar to those which play a key role in some of the proofs in that paper.

Let $A = \{a_1, \cdots a_m\}$ and $B = \{b_1, \cdots, b_n\}$, with the $a$'s and $b$'s as in Theorem 4. Again we assume w.l.o.g. that, in the statement of Theorem 4, the circuit $C$ bounds the outer face in the given embedding of $G$, that the vertices of $C$ are given in clockwise order, and that $U = \{u\}$, $W = \{w\}$. In this proof we will not use the notion of $\overline{G}$ and that of a drawing of $G$. In the description below we always have in mind the embedding of $G$ given in the statement of the theorem, with the above-mentioned assumptions.

Each path $\pi$ from $u$ to $w$ partitions the set $E(G)$ into three subsets: the set of edges of $\pi$ itself, the edges in the part of $G$ to the left of $\pi$ (when we follow $\pi$ in the direction of $w$), and the edges in the part of $G$ to the right of $\pi$. We denote these sets by $E(\pi)$, $E_L(\pi)$ and $E_R(\pi)$ respectively.

For each configuration $\omega \in \{0, 1\}^{E(G)}$ which has an open path from $u$ to $w$, we will consider the left-most self-avoiding path from $u$ to $w$ (this is similar to the well-known notion of lowest crossing in, e.g., bond percolation on a box in the square lattice). Analogously we will consider the right-most self-avoiding open path. For brevity we will drop the word 'self-avoiding'.

Let $\mathcal{P}$ denote the measure on $\{0, 1\}^{E(G)}$ corresponding to the given bond percolation model; that is, the product measure with parameters $p_e, e \in E(G)$. Let $\pi$ be a path from $u$ to $w$.

**Observation.** Conditioned on the event that $\pi$ is the leftmost open path from $u$ to $w$, each edge $e$ in $E_R(\pi)$ is, independently of the others, open with probability $p_e$.

An analogous observation holds when we replace 'leftmost' by 'rightmost' and $E_R(\pi)$ by $E_L(\pi)$.

It is easy to check that similar properties hold for the distribution $\mu$ obtained from $\mathcal{P}$ by conditioning on having an open path from $u$ to $w$:

$$\mu(\omega_e = \cdot, \, e \in E_R(\pi) \, | \, \pi \text{ is the leftmost open path from } u \text{ to } w)$$

is the product distribution on $\{0, 1\}^{E_R(\pi)}$ with parameters $p_e, e \in E_R(\pi)$, and similarly if we replace leftmost by rightmost, and $E_R$ by $E_L$.

Let $\Gamma$ be the set of all configurations $\omega$ that have an open path from $u$ to $w$. We will construct a Markov chain $\omega_n, n = 0, 1, \cdots$, with state space $\Gamma$ and stationary distribution $\mu$. To do this, we first introduce auxiliary $0 - 1$ valued random variables $l_n(e), r_n(e)$, $e \in E(G)$, $n = 0, 1, \cdots$. These random variables are independent and, for each $e$ and $n$,

$$P(l_n(e) = 1) = P(r_n(e) = 1) = p_e.$$

As initial state of the Markov chain we take $\omega_0 = \alpha$, for some $\alpha \in \Gamma$; the precise choice does not matter.

The transition from time $n$ to time $n + 1$ of this Markov chain consists of two substeps, (i) and (ii) below:

*Substep* (i): Denote by $\pi_n$ the leftmost open path from $u$ to $w$ in the configuration $\omega_n$. Using the $r_n$ variables introduced above, we update all edges in $E_R(\pi_n)$. This gives a new configuration, which we denote by $\omega'_n$. More precisely, we define

$$\omega'_n(e) = r_n(e), \quad \text{if } e \in E_R(\pi_n),$$

and $\omega_n(e)$ otherwise.

*Substep* (ii): To $\omega'_n$ we apply, informally speaking, the same action as in substep (i), but now with 'left' and 'right' exchanged. The resulting configuration is $\omega_{n+1}$. More precisely, with $\pi'_n$ denoting the rightmost open path from $u$ to $w$ in $\omega'_n$, we define

$$\omega_{n+1}(e) = l_n(e), \quad \text{if } e \in E_L(\pi'_n),$$

and $\omega'_n(e)$ otherwise.

Let $\mu_n$ denote the distribution of $\omega_n$. From the above construction it is clear that $\omega_n$, $n = 0, 1, \cdots$ is indeed a Markov chain with state space $\Gamma$. Moreover, it is clear from the construction and the above-mentioned Observation that $\mu$ is invariant under the above dynamics. (In fact, it is invariant under substep (i) as well as under substep (ii)). It is also easy to see that this Markov chain is aperiodic and irreducible. Hence, $\mu_n$ converges to $\mu$. Let, for each vertex $x$, $\eta_n(x)$ be the indicator of the event that the configuration $\omega_n$ has an open path from $u$ to $x$.

By the above arguments, it is sufficient to show that, for each $n$, the collection of random variables

$$\{\eta_n(x), x \in A\} \cup \{1 - \eta_n(y), y \in B\}$$

is positively associated. This, in turn, follows from the following Claim and the well-known Harris-FKG theorem that independent random variables are positively associated:

**Claim.** *Fix the initial configuration $\alpha$ of the Markov chain. For each $n$ and for each $x \in A$, $\eta_n(x)$ is then a function of the variables $l_k(e)$ and $r_k(e)$, $e \in E(G)$, $0 \le k \le n - 1$. Moreover, it is increasing in the $l$ variables and decreasing in the $r$ variables. An analogous statement, but with $l$ and $r$ interchanged, holds for $x \in B$.*

We give a brief sketch of the proof of this claim: Consider the following partial order, called 'more leftish than', on $\Gamma$. First let, for $\omega \in \Gamma$, $\pi_L(\omega)$ and $\pi_R(\omega)$ denote the leftmost and the rightmost open paths from $u$ to $w$ respectively. If $\omega$, $\hat{\omega} \in \Gamma$, we say that $\hat{\omega}$ is more leftish than $\omega$ iff each of the following ((a) and (b)) holds:

(a) The leftmost and the rightmost open path of $\hat{\omega}$ are located to the left of the corresponding paths of $\omega$. More precisely,

$$E(\pi_L(\hat{\omega})) \subset E_L(\pi_L(\omega)) \cup E(\pi_L(\omega)),$$

and

$$E(\pi_R(\omega)) \subset E_R(\pi_R(\hat{\omega})) \cup E(\pi_R(\hat{\omega})).$$

(b)

$$\hat{\omega}(e) \ge \omega(e), e \in E_L(\pi_L(\hat{\omega})),$$

and

$$\hat{\omega}(e) \le \omega(e), e \in E_R(\pi_R(\omega)).$$

With fixed $l$'s and $r$'s this order is preserved under substep (i) as well as under substep (ii). Moreover, it is easy to check that if we apply substep (i) or substep (ii) to some configuration $\omega$, and increase some of the $l$ variables or decrease some of the $r$ variables involved in that substep, the configuration resulting from that substep will become more leftish. These arguments, together with the fact that if $\hat{\omega}$ is more leftish than $\omega$, then

$$I_{\{u \to x \text{ in } \hat{\omega}\}} \ge I_{\{u \to x \text{ in } \omega\}}, \, x \in A,$$

and

$$I_{\{u \to x \text{ in } \hat{\omega}\}} \leq I_{\{u \to x \text{ in } \omega\}}, \, x \in B,$$

imply the above Claim.

## 4. Proof of Theorem 2

Consider the contact process in the statement of Theorem 2. A useful and well-known way to describe this process is by a so-called space-time diagram, or graphical representation (see e.g. figure 1 in part I of [7]): We represent each site $x$ as the point $(x, 0)$ in the plane and we assign to it a vertical line (time axis) $l_x = \{(x, t) \, : \, t > 0\}$. On $l_x$ we consider three independent Poisson point processes: one with density $\delta_x$, corresponding to recovery attempts; one with density $\lambda_{x-1,x}$ corresponding to attempts to infect site $x-1$; and one with density $\lambda_{x+1,x}$ corresponding to attempts to infect site $x+1$. At each point in the first point process we draw a symbol $*$ on $l_x$; from each point in the second point process we draw a horizontal arrow to $l_{x-1}$; similarly, from each point in the third point process we draw a horizontal arrow to $l_{x+1}$. By an *allowable* path we mean a continuous trajectory along the $l_x$'s and the arrows specified above, which satisfies the following conditions: along the $l_x$'s it goes only upward, and is not allowed to cross a $*$; when it goes along an arrow it must respect the direction of that arrow.

A site $y$ is infected at time $t$ iff for some site $x$ that is infected at time 0, there is an allowable path from $(x, 0)$ to $(y, t)$ in the space-time diagram; that is,

$$\eta_t(y) = I_{\{\exists x \text{ s.t. } x \text{ is infected at time 0 and}}$$
$$\text{there is an allowable path from } (x,0) \text{ to } (y,t)\}.$$

We will apply Theorem 4 to a discrete-time approximation of this process. Similar discretization arguments for contact processes (and many other interacting particle systems) are quite common (see e.g. pages 11 and 65 of [7]). Let $N$ be a positive integer. Consider bond percolation on the following graph, $G$. The vertex set of $G$ is $\{(x, k/N) \, : \, x \in \mathbf{Z}, k = 0, 1, \cdots\}$. Each vertex $(x, k/N)$ is the starting point of three oriented edges: one to $(x+1, k/N)$, one to $(x-1, k/N)$, and one to $(x, (k+1)/N)$. We take these edges to be open with probabilities $\lambda(x+1,x)/N$, $\lambda(x-1,x)/N$ and $1 - \delta_x/N$ respectively. Let $U = \{(u, 0) \, : \, u \text{ is infected at time 0}\}$. Now fix a $t > 0$ and a positive integer $n$. Let $\hat{t}$ be the smallest multiple of $1/N$ that is larger than or equal to $t$. Let $G_{n,N}$ be the (finite) subgraph induced by $V_{n,N} = \{(x, j/N) : |x| \leq n, j \leq N\hat{t}\}$.

Let, for each integer $x$, $\eta_t^{(n,N)}(x)$ denote the indicator of the event $\{U \to (x, \hat{t})$ in $G_{n,N}\}$. (It would be more correct here to write $U \cap ([-n, n] \times 0)$ in place of $U$.) Fix a positive integer $m$. It is quite standard that the joint distribution of the random variables $\eta_t^{(n,N)}(x)$, $-m \leq x \leq m$, converges to that of $\eta_t(x)$, $-m \leq x \leq m$, when we let $N \to \infty$ and then $n \to \infty$.

Moreover, to each graph $G_{n,N}$ we can apply Theorem 4, which tells us that, conditioned on the event $\{U \to (0, \hat{t})$ in $G_{n,N}\}$, the collection of random variables

$$\{1 - \eta_t^{(n,N)}(x) \, : \, -m \leq x \leq -1\} \cup \{\eta_t^{(n,N)}(x) \, : \, 1 \leq x \leq m\}$$

is positively associated. This, combined with the above-mentioned limit considerations, gives us that, conditioned on the event that $\eta_t(0) = 1$, the collection

$$\{1 - \eta_t(x) \, : \, -m \leq x \leq -1\} \cup \{\eta_t(x) \, : \, 1 \leq x \leq m\}$$

is positively associated. Since this holds for all $m$, Theorem 2 follows. $\square$

**Acknowledgment**

We thank N. Konno for bringing his conjecture to our attention.